\documentclass[leqno,12pt]{article}

\begin{document}

\title{Report on the trimestre ``Heat Kernels, Random Walks, and
Analysis on Manifolds and Graphs'' at the Centre \'Emile Borel 
(Institut Henri Poincar\'e, Spring, 2002)}

\author{Stephen Semmes}

\date{}

\maketitle


\renewcommand{\thefootnote}{}   

\footnotetext{Some historical notes, mentioned by a colleague: \'Emile
Borel spoke at the opening of the author's home institution, Rice
University (originally the Rice Institute) in Houston, Texas, in
1912.  Borel published ``Molecular theories and mathematics'' in
connection with his lectures in the Rice Institute Pamphlet, Volume I
(1915), 163--193.  Henri Poincar\'e was also invited by President
Edgar Odell Lovett and accepted, conditioned on the state of his
health, but eventually declined the invitation and subsequently passed
away.  Borel's paper begins with a tribute to Poincar\'e, and relates
a discussion they had about the trip.  Borel indicates that he would
have changed his subject to an appreciation of Poincar\'e's work,
except that Vito Volterra was doing exactly that.  Volterra's paper
appears in the same issue of the Rice Institute Pamphlet, ``Henri
Poincar\'e'', pp. 133--162.  Jacques Hadamard contributed ``The early
scientific work of Henri Poincar\'e'' and ``The later scientific work
of Henri Poincar\'e'' to the Rice Institute Pamphlet, Volume IX
(1922), 111-183 and Volume XX (1933), 1--86.  Hadamard makes the point
in the introduction to the first paper that uses for Poincar\'e's work
seemed to take 25 years to be found.}

\subsubsection*{``If it's Tuesday, this must be Belgium.''}  

This is the name of a film about a group of tourists who were going
from city to city a little too fast.  Fortunately in this trimestre
there was more time, and while the activities were numerous and
extensive, one had the opportunity to delve into various topics in
some detail.

	To give a part of the mathematical setting, let us review a
few classical matters related to calculus and partial differential
equations.  Fix a positive integer $n$, and let ${\bf R}^n$ be the
usual $n$-dimensional Euclidean space, consisting of $n$-tuples of
real numbers.  If $f(x)$ is a real-valued function on ${\bf R}^n$
which is twice-continuously differentiable, say, then the
\emph{Laplacian} of $f$ is denoted $\Delta f$ and defined by
\begin{equation}
	\Delta f = \sum_{j=1}^n \frac{\partial^2}{\partial x_j^2} f.
\end{equation}

	Let $f_1(x)$, $f_2(x)$ be two real-valued functions on ${\bf
R}^n$ which are continuous and have compact support, so that they are
both equal to $0$ outside of a bounded set.  More generally, one can
assume that $f_1$, $f_2$ satisfy suitable decay conditions, etc.  The
standard inner product of such functions is defined by
\begin{equation}
\label{inner product of functions}
	\langle f_1, f_2 \rangle = \int_{{\bf R}^n} f_1(x) \, f_2(x) \, dx.
\end{equation}
There is another symmetric bilinear form which is closely related to
the Laplacian, given by
\begin{equation}
\label{def of mathcal{E}(f_1,f_2)}
	\mathcal{E}(f_1,f_2) = 
	\frac{1}{2} \int_{{\bf R}^n} \nabla f_1(x) \cdot \nabla f_2(x) \, dx
\end{equation}
when $f_1$, $f_2$ are continuously differentiable, or satisfy other
appropriate regularity conditions.  Here $\nabla f(x)$ denotes the
gradient of $f$ at $x$, i.e., the vector with components $(\partial /
\partial x_j) f(x)$, and $v \cdot w$ is the usual inner product on
${\bf R}^n$, so that $v \cdot w = \sum_{j=1}^n v_j \, w_j$.  If in
addition $f_1$ is twice continuously-differentiable, then
\begin{equation}
	\mathcal{E}(f_1,f_2) = 
	- \frac{1}{2} \int_{{\bf R}^n} \Delta f_1(x) \, f_2(x) \, dx.
\end{equation}
This follows from integration by parts.

	The \emph{energy} $\mathcal{E}(f)$ of a function $f$ is
defined by
\begin{equation}
\label{def of mathcal{E}(f)}
	\mathcal{E}(f) = \mathcal{E}(f,f) =
	\frac{1}{2} \int_{{\bf R}^n} |\nabla f(x)|^2 \, dx,
\end{equation}
where $|v|$ denotes the standard Euclidean length of $v$, which is
the same as saying that $|v|^2 = v \cdot v$.  If $\eta(x)$ is another
function on ${\bf R}^n$, then
\begin{equation}
	\frac{d}{ds} \mathcal{E}(f + s \, \eta) \Big\vert_{s=0} =
		- \int_{{\bf R}^n} \Delta f(x) \, \eta(x) \, dx,
\end{equation}
under suitable conditions on $f$ and $\eta$.  This is commonly
rephrased as saying that the gradient of the enrgy functional
$\mathcal{E}(f)$ is given by $-\Delta f$, where this statement
implicitly uses the inner product (\ref{inner product of functions})
on functions on ${\bf R}^n$.

	A function $u(x,t)$ on ${\bf R}^n \times (0,\infty)$ which is
twice-continuously differentiable in $x$ and continuously
differentiable in $x$ and $t$ is said to satisfy the \emph{heat
equation} if
\begin{equation}
	\frac{\partial}{\partial t} u = \Delta u.
\end{equation}
Under modest growth conditions on a function $f(x)$ on ${\bf R}^n$,
there is a unique continuous function $u(x,t)$ on ${\bf R}^n \times
[0,\infty)$ such that $u(x,0) = f(x)$, $u(x,t)$ is infinitely
differentiable in $x$ and $t$ when $t > 0$, $u(x,t)$ satisfies the
heat equation on ${\bf R}^n \times (0,\infty)$, and $u(x,t)$ also
satisfies modest growth conditions (which can be related to those of
$f$).

	One way to look at the heat equation is as an ordinary
differential equation in $t$, acting in vector spaces of functions
of $x$.  To find $u(x,t)$ given $f(x)$ as in the preceding
paragraph, one might write
\begin{equation}
\label{u(x,t) = (exp (t Delta)f)(x)}
	u(x,t) = (\exp (t\Delta)f)(x).
\end{equation}
In fact the Fourier transform gives a useful way to make sense of this.

\subsubsection*{Aspects of symmetry}

	Versions of these notions come up in a variety of situations,
and a number of these were discussed in the trimestre.  In the spirit
of the book ``Introduction to Fourier Analysis on Euclidean Spaces''
by E.~Stein and G.~Weiss, which also provides a lot of helpful
background information for these topics, one might start by
considering the symmetries of the objects just described.  They are
all invariant under translations, and under rotations on ${\bf R}^n$.
They also behave nicely with respect to dilations on ${\bf R}^n$,
which is to say under transformations of the form $x \mapsto a x$,
where $a$ is a positive real number.  In the case of the heat
equation, one should use the dilations $(x,t) \mapsto (a x, a^2 t)$,
to adjust for the fact that there is one derivative in $t$ and
derivatives of order $2$ in $x$.

	Instead of Euclidean spaces a basic setting is that of
irreducible symmetric spaces of noncompact type, which was discussed
in the course of J.-P.~Anker.  For these one again has translation
invariance and forms of rotation invariance, but no dilation
invariance.  There are counterparts of Fourier analysis here too,
for analyzing solutions to the heat equation, but this has some
weaknesses differing from the Euclidean case.  

	In the Euclidean case the solution $u(x,t)$ to the heat
equation with initial data $f(x)$ can be expressed in the form
\begin{equation}
\label{u(x,t) = int_{{bf R}^n} k_t(x-y) f(y) dy}
	u(x,t) = \int_{{\bf R}^n} k_t(x-y) \, f(y) \, dy
\end{equation}
for a function $k_t(x)$ called the \emph{heat kernel}.  The fact
that the solution can be written in this manner, instead of
\begin{equation}
\label{u(x,t) = int_{{bf R}^n} k_t(x,y) f(y) dy}
	u(x,t) = \int_{{\bf R}^n} k_t(x,y) \, f(y) \, dy,
\end{equation}
reflects the translation-invariance of the problem in $x$.  The
rotation-invariance of the problem implies in turn that $k_t(x)$ is a
radial function of $x$, so that $k_t(x)$ can be written as $h_t(|x|)$
for a function $h_t(r)$ with $t \in (0,\infty)$ and $r \in
[0,\infty)$.  One can go further and use dilation-invariance to obtain
that $k_t(x)$ is of the form $t^{-n/2} h(|x|/\sqrt{t})$ for a function
$h(r)$, $r \in [0,\infty)$.  It is a classical result, which is a good
exercise to derive, that $k_t(x)$ is in fact a Gaussian function of
$x$.  This can be viewed in terms of the Fourier transform, or by
working out an ordinary differential equation for the function $h(r)$.

	In the context of symmetric spaces one can start with a
general form for $u(x,t)$ as in (\ref{u(x,t) = int_{{bf R}^n} k_t(x,y)
f(y) dy}), and use translation-invariance to reduce to something more
like (\ref{u(x,t) = int_{{bf R}^n} k_t(x-y) f(y) dy}).  The
counterpart of rotation-invariance permits one to reduce the number of
variables further, but not in general to $2$ variables.  Fourier
analysis leads to interesting representations for the heat kernel,
but fundamental features concerning size and localization are not
always so clear from this representation.

	Now let us go in a different direction and suppose that we are
working on ${\bf R}^n$ again, but with a differential operator $L$
with variable coefficients in place of the Laplacian.  Specifically,
we assume that $L$ is of the form
\begin{equation}
	L = \sum_{j,m = 1}^n \frac{\partial}{\partial x_j} \, a_{j,m}(x)
				\, \frac{\partial}{\partial x_m},
\end{equation}
where $a_{j,m}(x)$ are bounded real-valued functions which satisfy
\begin{equation}
	a_{j,m}(x) = a_{m,j}(x)
\end{equation}
and
\begin{equation}
\label{|v|^2 le sum_{j,m = 1}^n a_{j,m}(x) v_j v_m}
	|v|^2 \le \sum_{j,m = 1}^n a_{j,m}(x) \, v_j \, v_m
\end{equation}
for all $v \in {\bf R}^n$.  In other words, $(a_{j,m}(x))_{j,m}$ are
positive-definite real symmetric matrices which are uniformly bounded
in $x$ and bounded from below in the sense of matrices by the identity
matrix.  Because the coefficients are allowed to depend on $x$, we
lose in general the invariance under translations, rotations, or
dilations, and the heat kernel should be written as $k_t(x,y)$, with
$x, y \in {\bf R}^n$ and $t > 0$, as in (\ref{u(x,t) = int_{{bf R}^n}
k_t(x,y) f(y) dy}).  However, there are vestiges of these invariances,
in that translations and rotations of $L$ lead to operators of the
same type, and similarly for dilations if one includes suitable
scale-factors.  While the precise form of the heat kernel may not be
easy to describe, one can try to show that it has many properties in
common with the Gaussian kernels in the case of the standard
Laplacian.

	One can go further and consider coefficients $a_{j,m}(x)$
which are not symmetric in $j$ or $m$, and perhaps not even
real-valued.  For the latter one can adjust (\ref{|v|^2 le sum_{j,m =
1}^n a_{j,m}(x) v_j v_m}) by taking the real part of the right side,
so that one still has ``uniform ellipticity''.  More generally one
can allow operators of order larger than $2$, and vector-valued functions
and systems of differential equations.  Questions related to these
situations were discussed in the courses of p.~Auscher and P.~Tchamitchian,
and of S.~Hofmann and A.~McIntosh.

	Note that it still makes sense to talk about
\begin{equation}
\label{exp (t L)}
	\exp (t L)
\end{equation}
in this type of situation, using spectral theory.  This works more
nicely when the coefficients $a_{j,m}(x)$ are real and symmetric, so
that the operator $L$ is self-adjoint (with a suitable choice of
domain).  Even without these conditions, one can define (\ref{exp (t
L)}), using resolvent integrals.  For that matter, one can define more
general functions of $L$, and part of the interest of the heat kernels
is that the exponentials (\ref{exp (t L)}) and related operators can
make good building blocks for studying other functions of $L$.

	On a connected Lie group $H$ one can again look at
second-order elliptic differential operators $L$ which are invariant
under translations, but in general $H$ can be noncommutative and one
should be careful to specify whether $L$ is invariant under left
translations, right translations, or both.  In the case of Lie groups
which are nilpotent, such as the Heisenberg groups, dilations can be
used in much the same manner as on Euclidean spaces to have an extra
degree of symmetry.  In the course of W.~Hebisch, solvable Lie groups
and operators on them were treated, for which there is a delicate
interplay between exponential growth on the one hand and having a fair
amount of commutativity around on the other hand.  

	S.~Lang gave a series of lectures concerning deep questions of
expansions for heat kernels on the locally symmetric spaces (of finite
volume)
\begin{equation}
	SL(n,{\bf R}) / SL(n,{\bf Z}), \quad
		SL(n,{\bf C}) / SL(n,{\bf Z}[i]),
\end{equation}
where ${\bf Z}$ denotes the set of integers, and ${\bf Z}[i]$
is the set of complex numbers whose real and imaginary parts are
integers.

\subsubsection*{Discrete settings}

	Let us consider ${\bf Z}^n$ now instead of ${\bf R}^n$.  If
$x$, $y$ are elements of ${\bf Z}^n$, let us say that $x$ and $y$ are
\emph{adjacent} if $|x-y| = 1$.  Thus $x$ and $y$ are adjacent if they
agree in all but one component, where they differ by $\pm 1$.  If
$f(x)$ is a function on ${\bf Z}^n$, define $A(f)$ on ${\bf Z}^n$ by
\begin{equation}
	A(f)(x) = \frac{1}{2n} \sum_{y \in {\bf Z}^n \atop |x-y| = 1} f(y),
\end{equation}
so that $A(f)(x)$ is the average of $f$ over the $2n$ elements of
${\bf Z}^n$ adjacent to $x$.

	The linear operator $A - I$ on functions on ${\bf Z}^n$, where
$I$ denotes the identity operator, is a discrete version of the
Laplacian.  This makes more sense if one writes the classical Laplacian
of a twice continuously-differentiable function $h$ at a point $x$ as
\begin{equation}
	\Delta (h)(x) = \lim_{r \to 0} \frac{1}{r^2} ({\rm Av}(h)(x,r) - h(x)),
\end{equation}
with ${\rm Av}(h)(x,r)$ equal to the average of $h$ over the sphere
with center $x$ and radius $r$.

	The analogue of the heat equation for a function $u(x,t)$ with
$x$ in ${\bf Z}^n$ and $t$ ranging through nonnegative integers can be
written as
\begin{equation}
	u(x,t+1) = \frac{1}{2n} \sum_{y \in {\bf Z}^n \atop |x-y| = 1} u(x,t),
\end{equation}
which is the same as saying that $u(x,t+1)$ is given by applying the
operator $A$ to $u(x,t)$ as a function of $x$.  To make this look more
like the classical heat equation, one can reexpress this as saying
that $u(x,t+1) - u(x,t)$, which is like the ``derivative'' of $u$ in
$t$, is equal to $A - I$ applied to $u(x,t)$ as a function of $x$.
Clearly, for any function $f(x)$ on ${\bf Z}^n$, there is a unique
function $u(x,t)$ defined for $x$ in ${\bf Z}^n$ and $t$ a nonnegative
integer such that $u(x,0) = f(x)$ for all $x$ in ${\bf Z}^n$ and $u(x,t)$
satisfies the heat equation above for all $x$ and $t$.  In fact,
$u(x,t)$ can be written as
\begin{equation}
	u(x,t) = (A^t)(f)(x),
\end{equation}
in analogy with (\ref{u(x,t) = (exp (t Delta)f)(x)}).

	In analogy with (\ref{u(x,t) = int_{{bf R}^n} k_t(x-y) f(y) dy}),
we can write
\begin{equation}
\label{u(x,t) = sum_{y in {bf Z}^n} p_t(x-y) f(y)}
	u(x,t) = \sum_{y \in {\bf Z}^n} p_t(x-y) \, f(y),
\end{equation}
where the ``heat kernel'' $p_t(w)$ is defined for $t$ a nonnegative
integer and $w$ in ${\bf Z}^n$.  Specifically, $p_0(w)$ is equal to
$0$ when $w \ne 0$ and to $1$ when $w = 0$, $p_1(w)$ is equal to
$0$ when $w$ is not adjacent to $0$ and to $1/(2n)$ when $w$ is
adjacent to $0$, and $p_t(w)$ can easily be determined explicitly.

	In fact, $p_t(x-y)$ is the probability that the standard
random walk on ${\bf Z}^n$ goes from $x$ to $y$ in exactly $n$ steps.
In the continuous setting there are similar statements for Brownian
motion and other processes associated to second-order differential
operators.

	That the heat kernel in (\ref{u(x,t) = sum_{y in {bf Z}^n}
p_t(x-y) f(y)}) is of the form $p_t(x-y)$, rather than $p_t(x,y)$,
reflects the translation-invariance here, just as in the classical
case on ${\bf R}^n$.  Of course one can consider other graphs instead
of ${\bf Z}^n$, with similar objects as defined above, and with a
formula of the type
\begin{equation}
	u(x,t) = \sum p_t(x,y) \, f(y)
\end{equation}
in place of (\ref{u(x,t) = sum_{y in {bf Z}^n} p_t(x-y) f(y)}).

	The course of T.~Sunada dealt with \emph{crystal lattices}, which
are characterized in terms of a large abelian group of symmetries.
The graphs ${\bf Z}^n$ are a very special case of this, and numerous
other configurations are possible.  In W.~Woess' course, techniques of
\emph{generating functions} were discussed, which can lead to remarkable 
formulas and information about random walks.  Part of M.~Barlow's course
was concerned with random walks on graphs with self-similarity, and
the effect of self-similarity on the heat kernel.

	In analogy with second-order differential operators on ${\bf
R}^n$ with variable coefficients, one can consider random walks and
discrete Laplacians on ${\bf Z}^n$ in which the weighting factors vary
from point to point.  One does not need to stick to ${\bf R}^n$ or
${\bf Z}^n$ here; one can work on manifolds or graphs, or more
generally metric spaces equipped with a measure.  Several of the courses
dealt with different facets of this, including Sobolev spaces and
Sobolev or Poincar\'e inequalities.

	R.~Brooks discussed in his course Riemann surfaces, graphs,
correspondences between them, and lower bounds for positive eigenvalues
for the Laplacian for both.

\subsubsection*{Additional topics}

	Let $p$ be a real number, $p > 1$.  For suitable functions
$f(x)$ on ${\bf R}^n$, consider the \emph{$p$-energy functional}
\begin{equation}
	\mathcal{E}_p(f) = \frac{1}{p}\int_{{\bf R}^n} |\nabla f(x)|^p.
\end{equation}
This is the same as $\mathcal{E}(f)$ in (\ref{def of mathcal{E}(f)})
when $p = 2$, but there is not a bilinear version as in (\ref{def of
mathcal{E}(f_1,f_2)}) when $p \ne 2$.  However, one can again consider
the derivative of $\mathcal{E}_p(f)$ in $f$ for all $p$, and this
leads to a nonlinear (when $p \ne 2$) second-order differential
operator known as the \emph{$p$-Laplacian}.

	The $p$-energy is invariant under translations and rotations,
and scales under dilations in a simple way, just as when $p = 2$.  For
$p = n$ there is additional symmetry, known as \emph{conformal
invariance}.

	One can consider more complicated functionals which behave in
roughly the same manner in terms of size, but which incorporate
``variable coefficients'' into the picture.  When $p = n$ there is a
``quasi-invariance'' of the energy under \emph{quasiregular} mappings,
which are defined in terms of a pointwise quasiconformality property
(where the $n$th power of the norm of the differential of the mapping
is bounded by a constant times the Jacobian, i.e., the determinant of
the differential of the mapping).  Quasiregular mappings, unlike
quasiconformal mappings, are allowed to have branching, analogous to
holomorphic mappings in the complex plane which are not one-to-one.
The \emph{quasi-invariance} of the $p$-energy when $p = n$ states that
the energy functional is transformed by a quasiregular change of
variables into an energy functional of roughly the same type, but with
variable coefficients which satisfy bounds in terms of the
quasiregularity constant.  As a result, a solution of the $n$-Laplace
equation is transformed, after composition with a quasiregular mapping,
into a solution of an analogous equation with variable coefficients,
still with suitable boundedness and ellipticity conditions.  This is
an important tool in the study of quasiregular mappings, as discussed
in the course of I.~Holopainen.

	Even with the extra nonlinearity, there are similar issues
concerning the relationship between the geometry of a space and the
behavior of solutions of differential equations or inequalities as
before.

	A different kind of nonlinearity was treated in the course of
K.-T.~Sturm, with averages, heat flows, and random processes taking
values in a metric space, under general conditions of nonpositive
curvature.  It can be clear how to take a weighted average of two
points in a metric space, using a point along a geodesic arc that
joins them, but for more than two points not lying on the same
geodesic the situation becomes more complicated.  A fascinating
feature of the probabilistic point of view is that in a sequence of
independent samples one can use the ordering of the sequence to apply
the two-point case step-by-step; it turns out that there are results
to the effect that the limit of this exists and is the same almost
surely, and that the common answer is the same as one produced from
another procedure which deals with all points in the average at the
same time.

	The courses of B.~Driver and L.~Saloff-Coste were concerned with
analysis on infinite-dimensional spaces.  Specifically, Driver's course
dealt with Weiner space, spaces of paths in manifolds, and loop groups,
while Saloff-Coste's course addressed locally-compact and connected
topological groups, such as infinite products of finite-dimensional
compact connected Lie groups.

	Of course the brief overview given here is not at all intended
to be exhaustive.  Fortunately, a volume is in preparation containing
surveys and other material from the trimestre, in which much more
information can be found.

\end{document}